\documentclass[leqno,12pt]{article}
\usepackage{amssymb,amsbsy,amsmath,amsfonts,amssymb,amscd}
\usepackage{latexsym}
\usepackage{times}
\begin{document}
\newtheorem{theorem}{Theorem}
\newtheorem{lemma}{Lemma}
\newtheorem{remark}{Remark}
\newtheorem{conjecture}{Conjecture}
\title{\Large\sc Green-Lazarsfeld's Conjecture for Generic Curves
of Large Gonality}
\author{\normalsize Marian Aprodu\footnote{Supported by an E.C. Marie Curie
Fellowship, contract number HPMF-CT-2000-00895}\hskip 3mm
and Claire Voisin}
\maketitle

\begin{abstract}
\noindent
We use Green's canonical syzygy conjecture for generic curves
to prove that the Green-Lazarsfeld
gonality conjecture holds for generic curves of genus $g$, and gonality
$d$, if $g/3<d<\left[g/2\right]+2$. 
\end{abstract}

\begin{center}
{\bf 1. Introduction}
\end{center}

\noindent
Denoting by $K_{p,q}(X,L)$ the Koszul cohomology with value in a line
bundle $L$ ({\em see} \cite{Gr}), Green and Lazarfeld proved the
following ({\em cf.} Appendix to \cite{Gr}):
\begin{theorem}
\label{theorem: GL}
Let $X$ be a complex manifold, $L_1$, and $L_2$ be two line
bundles on $X$ such that $r_1:=h^0(L_1)-1\geq 1$, and
$r_2:=h^0(L_2)-1\geq 1$. 
Then $K_{r_1+r_2-1,1}(X,L_1\otimes L_2)\not= 0$.
\end{theorem}

Let now $C$ be a smooth complex smooth projective curve of gonality 
$$
d:=\mbox{min}\{\mbox{deg}(L),h^0(L)\geq 2\}.
$$

Green-Lazarsfeld's theorem, applied to
$L_1={\cal O}_C(D)$, where $\mbox{deg}(D)=d$,\linebreak
$h^0(C,{\cal O}_C(D))=2$, and $L_2=L-D$ with 
$\mbox{deg}(L)$ sufficiently large, implies
$$
K_{h^0(L)-d-1,1}(C,L)\not= 0.
$$

The gonality conjecture predicts that this is optimal, namely:
\begin{conjecture}[{\rm Green and Lazarsfeld, \cite{GL}}]
\label{conjecture: GL}
For $C$ a curve of gonality $d$, and for any line bundle $L$ of
sufficiently large degree, we have  $K_{h^0(L)-d,1}(C,L)= 0$.
\end{conjecture}

In spite of the evidence brought by the ``$K_{p,1}$ Theorem'' of
\cite{Gr}, which solves the hyperelliptic case (among other things),
after having formulated the gonality conjecture, Green and Lazarsfeld
shown their mistrust of the statement they had just made. Since then,
the conjecture has been almost forgotten, and it took a while untill some
new evidence was discovered ({\em see} \cite{Ehb}, \cite{Ap}). 
This delay is probably due to the fact that the conjecture did not count
among the mathematical highlights of the last years, as almost all the
attention in the theory of syzygies of curves was focused on the more
famous Green conjecture.

The aim of this short note is to mix together the main results of
\cite{Ap} on the one hand, and of \cite{Vo}, and \cite{Vo2} on the
other hand, in order to verify the Green-Lazarfeld conjecture for generic
curves of large given gonality. The first result we prove is the
following:
\begin{theorem}
\label{theorem: first}
For any positive integers $g$ and $d$ such that
$g/3+1\leq d<\left[(g+3)/2\right]$, the gonality conjecture is valid for
generic curves of genus $g$ and gonality $d$.
\end{theorem}

Note that for generic curves of genus $g$ the gonality equals $[(g+3)/2]$,
and thus Theorem \ref{theorem: first} covers all possible, not too small,
gonalities, except for the generic gonality. 

Our second result is:
\begin{theorem}
\label{theorem: second}
The gonality conjecture is valid for generic curves of 
even genus.
\end{theorem}

In the statements above, the word {\em generic} should be read in the
usual sense. The complex curves of fixed genus $g$, and gonality $d$, are
parametrised by an irreducible subvariety of the moduli space 
${\cal M}_g$, and a {\em generic curve} is a curve corresponding to a
generic point of this variety. Irreducibility follows from the
well-known fact that the closure of this subvariety is actually the
closure of the image in ${\cal M}_g$ of a Hurwitz scheme 
and then apply \cite{F}. 

Finally, we mention that all the notation we use in the sequel is
standard, and we refer to \cite{Gr} for basic facts about Koszul
cohomology.

\begin{center}
{\bf 2. Proofs of main results}
\end{center}

\noindent
First of all, we recall the following result from \cite{Ap}:
\begin{theorem}
\label{theorem: Ap}
If $L$ is a nonspecial line bundle on a curve 
$C$, which satisfies $K_{n,1}(C,L)=0$, for a positive integer
$n$, then, for any effective divisor $E$ of degree $e$,
we have\linebreak
$K_{n+e,1}(C,L+E)=0$.
\end{theorem}

In particular, if $K_{h^0(L)-d,1}(C,L)=0$ for a nonspecial line bundle
$L$, with $h^0(L)-d>0$, then $K_{h^0(L^\prime)-d,1}(C,L^\prime)=0$ for any
$L^\prime$ of sufficiently large degree. By means of the Zariski
semi-continuity of graded Betti numbers ({\em see}, for example, 
\cite{BG}), for both Theorem \ref{theorem: first}, and Theorem
\ref{theorem: second}, it suffices to exhibit one $d$-gonal curve $C$ of
genus $g$ (where $d=g/2+1$ for Theorem \ref{theorem: second}), and one
nonspecial line bundle
$L$ on $C$ satisfying $K_{h^0(L)-d,1}(C,L)=0.$
\\\\
{\em Proof of Theorem \ref{theorem: first}.}
{\em Step 1.} Construction of $C$ and $L$.

\bigskip

A suitable choice of such a curve $C$ is provided by the proof of Corollary 1 
of \cite{Vo}. We start with a $K3$ surface $S$ whose Picard group is 
cyclic, generated by a line bundle $\mathcal L$ of self-intersection
${\mathcal L}^2=4k-2$, where $k=g-d+1$. We denote $\nu=g-2d+2\geq 1$.
Under these assumptions, as $\nu\leq k/2$, we know that there
exists an irreducible curve $X\in |\mathcal L|$, having exactly $\nu$
simple nodes as singular points, and no other singularities, and such
that its normalization $C$ is of gonality $d=k+1-\nu$, {\em see}
\cite{Vo}. 

We set $L=K_C(p+q)$, where $p$, and $q$ are two distinct points of $C$
that lie over a node $x$ of $X$.
\\\\
{\em Step 2.} Recall the main result of \cite{Vo}.
\begin{theorem}
\label{theorem: Vo}
The $K3$ surface $S$ being as above, we have 
$K_{k,1}(S,{\cal L})=0$.
\end{theorem}

\bigskip

It follows directly from this result, from the adjunction formula, and
from the hyperplane section theorem for Koszul cohomology ({\em see}
\cite{Gr}) that $K_{k,1}(X,K_X)=0$. Since $k=h^0(C,L)-d$, the proof of
our theorem is concluded by the following lemma.
\begin{lemma}
\label{lemma: inclusion}
Let $X$ be a nodal curve, 
$C\stackrel{f}{\longrightarrow}X$ be the normalization 
of $X$, and $p,q\in C$ be two distinct points lying over the same node
$x$ of $X$. Then, for any $n\geq 1$, we have a natural injective map
$K_{n,1}(C,K_C(p+q)) \subset K_{n,1}(X,K_X)$.
\end{lemma}

\noindent
{\em Proof of Lemma \ref{lemma: inclusion}.}
Firstly, we remark that there is a natural inclusion of spaces of sections
$H^0(C,K_C(p+q)) \subset H^0(X,K_X)$. Indeed, the two spaces are both 
contained in the space of meromorphic differentials on $C$. Thus
$H^0(X,K_X)$ identifies to the meromorphic differentials on $C$
having poles of multiplicity one over the nodes, and regular outside
these points, and whose sums of residues over the nodes vanish.
The inclusion above is then a direct consequence of the Residue
Theorem.

This inclusion yields the following injection between the Koszul complexes
of
$K_C(p+q)$ and
$K_X$

$$
\begin{matrix}
0\rightarrow &
\bigwedge^{n+1}H^0(C,K_C(p+q))&{\rightarrow}
&\bigwedge^nH^0(C,K_C(p+q))\otimes H^0(C,K_C(p+q))&\rightarrow\ldots\\
&\downarrow&&\downarrow&\\
0\rightarrow &\bigwedge^{n+1}H^0(X,K_X)&{\rightarrow}
&\bigwedge^nH^0(X,K_X)\otimes H^0(X,K_X)&\rightarrow\ldots
\end{matrix}
$$

To conclude that this induces an injection on the degree $1$
cohomology groups, we use the existence of the retraction-homotopy
(up to a coefficient of $(n+1)!$)  given by the wedge product :

$$
\bigwedge^nH^0(C,K_C(p+q))\otimes H^0(C,K_C(p+q))\rightarrow 
\bigwedge^{n+1}H^0(C,K_C(p+q)).
$$

\noindent
{\em Proof of Theorem \ref{theorem: second}.} {\em Step 1.} Construction 
of $C$ and $L$.

\bigskip

We make use of the same curves as those used in \cite{Vo2}. Let $S$ be a
$K3$ surface whose Picard group is generated by an ample line bundle 
$\cal L$ of self-intersection ${\cal L}^2=4k$, where $g=2k$, and by a
rational curve $\Delta$ such that ${\cal L}.\Delta =2$. We
choose $X$ an irreducible nodal curve in the linear system $|\cal L|$,
having exactly one node as singularity. The curve $X$ has arithmetic
genus $2k+1$. We denote by $C$ the normalization of $X$, $p$ and $q$ the
two distinct points of $C$ lying over the node of $X$, and $L=K_C(p+q)$. 
Thus $C$ has genus $2k$.
\\\\
{\em Step 2.} Recall the following result from \cite{Vo2}:
\begin{theorem}
\label{theorem: Vo2}
With the notation above, 
we have $K_{k,1}(S,{\cal L})=0$.
\end{theorem}

Applying the hyperplane section theorem we conclude that
$K_{k,1}(X,K_X)=0$. By means of Lemma \ref{lemma: inclusion}, we obtain
$K_{k,1}(C,L)=0$. Since $k=h^0(L)-(k+1)$, it follows that the curve $C$
is of maximal gonality $k+1$, and that the gonality conjecture is valid
for $C$.

\begin{center}
{\bf 3. Final remarks}
\end{center}

\begin{remark}
{\rm Using Theorem \ref{theorem: Ap} one can see that for any of the
curves $C$ considered in the proofs of Theorem \ref{theorem: first}, and
Theorem \ref{theorem: second}, the vanishing
$K_{h^0(L^\prime)-d,1}(C,L^\prime)=0$ predicted by the gonality
conjecture holds for any line bundle $L^\prime$ of degree at least $3g$.
The same is true for generic $d$-gonal curves of genus $g$, where $g$ and
$d$ satisfy $g/3<d<\left[g/2\right]+2$.}
\end{remark}
\begin{remark}
{\rm The case $g>(d-1)(d-2)$ has already been treated in
\cite{Ap}. Therefore, from the viewpoint of the Green-Lazarfeld
conjecture for generic curves of fixed genus $g$, and gonality $d$, for
$d\geq 6$, there is a gap remaining for $3d-3<g\leq (d-1)(d-2)$ that has
to be solved differently. The case of maximal gonality
$d=k+2$ in the odd genus case $g=2k+1$, which seems to be the most
difficult, is also left over. Nevertheless, we should point out that in
all these excepted cases, gonality conjecture is {\em almost true}, that
is, any line bundle $L$ of degree at least $3g-2$ on a generic
$d$-gonal curve $C$ of genus $g$ (for any choice of $g$ and $d$) satisfies
$K_{h^0(L)-d+1,1}(C,L)=0$. This follows directly from the main results of
\cite{T}, \cite{Vo}, \cite{Vo2}, and from Theorem 3 of \cite{Ap}.}
\end{remark}
\begin{remark}
{\rm An alternative proof of Lemma \ref{lemma: inclusion} can be obtained
in a more  algebraic way, by factoring the normalization morphism through
$C\stackrel{h}{\rightarrow}Y\stackrel{g}{\rightarrow}X$,
where $g$ is the smoothing of the node $x$. We analyse the
three morphisms between the structure sheaves, and we obtain
an isomorphism $H^0(X,K_X)\cong H^0(Y,g^*K_X)$, and
an inclusion $H^0(C,K_C(p+q))\subset H^0(Y,g^*K_X)$. Denoting
by $W=H^0(C,K_C(p+q))$ inside $H^0(X,K_X)$, we have
$K_{n,1}(C,K_C(p+q))\subset K_{n,1}(X,K_X,W)$. Next, we
use the embedding $K_{n,1}(X,K_X,W)\subset K_{n,1}(X,K_X)$,
which arises from the spectral sequence of \cite{Gr2}, to conclude.

This shows that in Lemma \ref{lemma: inclusion} the fact of working over
the complex numbers is not essential, and the statement is true for nodal
curves over any algebraically closed field.}
\end{remark}
\begin{remark}
{\rm The use of the results of \cite{Vo}, \cite{Vo2} in the
proofs of our theorems shades a new light on the relationships between 
Green's conjecture, and the gonality conjecture, the two statements seemed
to be (in a somewhat misterious way) intimately related to each other
({\em see} also \cite{Ap}). In fact, Lemma \ref{lemma: inclusion} answers 
partially to the conjecture made in \cite{Ap}.}
\end{remark}

\noindent
{\sc Adresses}
\\\\
{Marian Aprodu 
\\
Romanian Academy, Institute of Mathematics
"Simion Stoilow", P.O.Box 1-764, RO-70700, 
Bucharest, Romania (e-mail:
Marian.Aprodu\char64 imar.ro)
\\
and Universit\'e de Grenoble 1,
Laboratoire de Math\'ematiques,
Institut Fourier BP 74,
F-38402 Saint Martin d'H\`eres Cedex,
France (e-mail: aprodu\char64 mozart.ujf-grenoble.fr)
\\\\
Claire Voisin
\\
Universit\'e Paris 7 Denis Diderot - CNRS UMR 7586 -
Institut de Math\'ematiques 2, Place Jussieu, F-75251
Paris Cedex 05 (e-mail: voisin\char64 math.jussieu.fr)

\end{document}